\documentclass[12pt,a4paper]{article}
\usepackage[utf8x]{inputenc}
\usepackage[english]{babel}
\usepackage[T1]{fontenc}
\usepackage{amsmath}
\usepackage{amsfonts}
\usepackage{amssymb}
\usepackage{makeidx}
\usepackage{graphicx}
\usepackage[left=3cm,right=2cm,top=3cm,bottom=2cm]{geometry}
\usepackage{setspace}
\usepackage{amsthm}
\usepackage{cite}

\usepackage[colorlinks=true,urlcolor=blue,
citecolor=red,linkcolor=blue,linktocpage,pdfpagelabels,
bookmarksnumbered,bookmarksopen]{hyperref}

\addto\captionsenglish{}

\newtheorem{definition}{Definition}[section]

\newtheorem{proposition}{Proposition}[section]

\newtheorem{theorem}{Theorem}[section]

\newtheorem{example}{Example}[section]

\newtheorem{lemma}{Lemma}[section]

\newtheorem{observation}{Remark}[section]

\newtheorem{corollary}{Corolary}[section]

\numberwithin{equation}{section}

\newcommand{\cupdot}{\mathbin{\mathaccent\cdot\cup}}

\newcommand{\bo}{\begin{observation}}
\newcommand{\eo}{\end{observation}}
\newcommand{\bd}{\begin{definition}}
\newcommand{\ed}{\end{definition}}
\newcommand{\bp}{\begin{proposition}}
\newcommand{\ep}{\end{proposition}}
\newcommand{\bt}{\begin{theorem}}
\newcommand{\et}{\end{theorem}}
\newcommand{\bc}{\begin{corollary}}
\newcommand{\ec}{\end{corollary}}
\newcommand{\bl}{\begin{lemma}}
\newcommand{\el}{\end{lemma}}
\newcommand{\be}{\begin{example}}
\newcommand{\ee}{\end{example}}
\newcommand{\beq}{\begin{equation}}
\newcommand{\eeq}{\end{equation}}
\newcommand{\beqa}{\begin{equation*}}
\newcommand{\eeqa}{\end{equation*}}
\newcommand{\R}{\mathbb{R}}
\newcommand{\RN}{\mathbb{R}^{N}}
\newcommand{\N}{\mathbb{N}}
\newcommand{\Rdois}{\mathbb{R}^{2}}
\newcommand{\Lp}{ L^{p}(\mathbb{R}) }
\newcommand{\Ldois}{ L^{2}(\mathbb{R}) }
\newcommand{\Lquatro}{L^{4}(\R)}

\newcommand{\A}{\mathcal{A}}

\newcommand{\Oo}{\mathcal{O}}
\newcommand{\Nn}{\mathcal{N}}

\newcommand{\Humeio}{H^{\frac{1}{2}}(\mathbb{R})}

\newcommand{\Lt}{ L^{t}(\mathbb{R})}

\newcommand{\intR}{\displaystyle\int\limits_{\mathbb{R}}}
\newcommand{\un}{u_{n}}

\newcommand{\unk}{u_{n_{k}}}

\newcommand{\uk}{u_{k}}

\newcommand{\vn}{v_{n}}

\newcommand{\yn}{y_{n}}
\newcommand{\wn}{w_{n}}

\newcommand{\until}{\tilde{u}_n}
\newcommand{\vntil}{\tilde{v}_n}
\newcommand{\RA}{\rightarrow}
\newcommand{\CF}{\rightharpoonup}
\newcommand{\Frac}{(-\Delta)^{\frac{1}{2}}}
\newcommand{\ds}{\displaystyle\int\limits}

\usepackage{times, color, xcolor}
\addto\captionsportuguese{
}

\begin{document}

 	\title{The  Choquard logarithmic 
equation involving fractional Laplacian operator and a  nonlinearity   with  exponential critical  growth 
		\thanks{The first author was supported    by  Coordination of Superior Level Staff Improvement-(CAPES) -Finance Code 001 and  S\~ao Paulo Research Foundation- (FAPESP), grant $\sharp $ 2019/22531-4,
while the second  author was supported by  National Council for Scientific and Technological Development -(CNPq),   grant $\sharp $ 307061/2018-3 and FAPESP  grant $\sharp $ 2019/24901-3.
			}}
	\author{
		Eduardo  de S. Böer \thanks{ E-mail address: eduardoboer04@gmail.com, Tel. +55.51.993673377}  and Ol\'{\i}mpio H. Miyagaki \thanks{Corresponding author} \footnote{ E-mail address: ohmiyagaki@gmail.com, Tel.: +55.16.33519178 (UFSCar).}\\
		{\footnotesize Department of Mathematics, Federal University of S\~ao Carlos,}\\
		{\footnotesize 13565-905 S\~ao Carlos, SP - Brazil}\\ }
\noindent
				
	\maketitle

\noindent \textbf{Abstract:} In the present work we investigate the existence and multiplicity of nontrivial solutions for the Choquard Logarithmic equation $(-\Delta)^{\frac{1}{2}} u + au + \lambda (\ln|\cdot|\ast |u|^{2})u = f(u) \textrm{ \ in \ } \mathbb{R}$, for $ a>0 $, $ \lambda >0 $ and a nonlinearity $f$ with exponential critical growth. We prove the existence of a nontrivial solution at the mountain pass level and a nontrivial ground state solution under  exponential critical and subcritical growth. Morever, when $ f $ has subcritical growth we guarantee the existence of infinitely many solutions, via genus theory.

\vspace{0.5 cm}


\section{Introduction}
In the present paper we are concerned with existence and multiplicity results for the Choquard logarithmic equation
\begin{equation} \label{P}
\Frac u + u + (\ln|\cdot|\ast |u|^{2})u = f(u) \textrm{ \ in \ } \mathbb{R},
\end{equation}
where $ a=\lambda=1 $, $ \Frac $ is the fractional laplacian, $ \ln $ is the neperian logarithm and $f: \R \RA \R $ is continuous with exponential growth and primitive $ F(s)=\int\limits_{0}^{s}f(t)dt $. 

We recall that a function $ h $ has \textit{subcritical} exponential growth at $ +\infty $, if
$$
\lim\limits_{s\RA + \infty}\dfrac{h(s)}{e^{\alpha s^2}-1} = 0 \textrm{ \ , for all \ } \alpha >0 ,
$$
and we say that $ h $ has \textit{critical} exponential growth at $ +\infty $, if there exists $ \omega \in (0, \pi] $ and  $ \alpha_0 \in (0, \omega) $ such that
$$
\lim\limits_{s\RA + \infty}\dfrac{h(s)}{e^{\alpha s^2}-1} = \left\{ \begin{array}{ll}
0, \ \ \ \forall \ \alpha > \alpha_{0}. \\
+\infty , \ \ \ \forall \ \alpha < \alpha_0 .
\end{array} \right. 
$$

So, in the following we present some necessary conditions to obtain our main results. This kind of hypothesis are usual in works with Moser-Trudinger inequality, such as \cite{[frac], [9]}. We assume that $ f $ satisfies
$$ f\in C(\R , \R), f(0)=0, \mbox{has critical exponential growth and} \ F(t) \geq 0 \mbox{ \ , for all \ } t\in \R . \leqno{(f_1)}$$
$$ \lim\limits_{|t|\RA 0} \dfrac{f(t)}{|t|}=0.  \leqno{(f_2)}$$

From $ (f_1) $ and $(f_2)$, for $ q>2 $,  $\varepsilon >0 $ and $\alpha > \alpha_0$, there exists a constant $ b_2 > 0 $ such that
\begin{equation}\label{eq26}
|F(u)| \leq \dfrac{\varepsilon}{p}|u|^2 + b_2 |u|^q (e^{\alpha |u|^2}-1)  \ , \ \forall \ u \in X .
\end{equation} 
and there exists a constant $ b_3 > 0 $ satisfying
\begin{equation}\label{eq27}
|f(u)| \leq \varepsilon |u|+b_3 |u|^{q-1}(e^{\alpha |u|^2}-1) \ , \ \forall \ u \in X .
\end{equation}
We also can derive from (\ref{eq26}) a very useful inequality. Consider $ r_1, r_2 >1 $, $ r_1 \sim 1 $ and $ r_2 > 2 $, such that $ \frac{1}{r_1}+\frac{1}{r_2}=1 $. Then, for all $ u\in \Humeio $, $ \alpha > \alpha_0 $ and $ \varepsilon >0 $, we get that
\begin{equation}\label{eq6}
\intR |F(u)| dx \leq \dfrac{\varepsilon}{2}||u||_{2}^{2}+b_2||u||_{r_2 q}^{q}\left(\intR (e^{r_1 \alpha u^2}-1) dx \right)^{\frac{1}{r_1}} .
\end{equation}

Our strategy to prove Theorem \ref{t1} will consist in finding a Cerami sequence for the mountain pass level. In order to verify that such sequence is bounded in $\Humeio$, we will need the following condition 
$$ \mbox{ there exists }\ \theta > 4 \ \mbox{ such that}\   f(t)t \geq \theta F(t) > 0, \ \mbox{for all }\  t\in \R . \leqno{(f_3)}$$

Moreover, as we will be working with an exponential term, to guarantee that the mentioned Cerami sequence and the minimizing sequence for the ground state satisfy the exponential estimative, we rely in the next condition.
$$ \mbox{ there  exist}\  q>4 \ \mbox{ and}\  C_q> \dfrac{[2(q-2)]^{\frac{q-2}{2}}}{q^{\frac{q}{2}}}\dfrac{S_{q}^{q}}{\rho_{0}^{q-2}} \ \mbox{ such that}\  F(t) \geq C_q |t|^q , \ \mbox{for all} \  t \in \R , \leqno{(f_4)}$$ 
for $ S_q, \rho_0 >0 $ to be defined in Lemma \ref{l10}.

Hence, we are able to enunciate our first main result.

\bt\label{t1}
Assume $ (f_1)-(f_4) $, $ q>4 $ and $ C_q>0 $ sufficiently large. Then,
\begin{itemize}
\item[(i)] Problem (\ref{P}) has a solution $ u\in X\setminus\{0\} $ such that $$ I(u)=c_{mp}=\inf\limits_{\gamma \in \Gamma}\max\limits_{t\in [0, 1]}I(\gamma(t)) ,  $$
where $ \Gamma = \{ \gamma \in C([0, 1], X) \ ; \ \gamma(0) \ , \ I(\gamma(1))< 0 \} $. 

\item[(ii)] Problem (\ref{P}) has a ground state solution $ u\in X\setminus\{0\} $ in the sense that $ I(u)=c_g = \inf\{ I(v) \ ; \ v\in X \mbox{ \ is a solution of (\ref{P})} \} $.
\end{itemize}
\et

For the second main result, we are concerned with multiplicity of solutions. However, to obtain this we need to exchange the condition $ (f_1) $ by the condition below.
$$ f\in C(\R , \R), \ f \mbox{ \ is odd, has subcritical exponential growth and} \ F(t) \geq 0 \mbox{, for all \ } t\in \R . \leqno{(f_1 ')}$$
Also, we need to add a condition that gives us the desired geometry for the associated functional. That is, 
$$ \mbox{the function \ } t\mapsto \dfrac{f(t)}{t^3} \mbox{ \ is increasing in \ } (0, + \infty) . \leqno{(f_5)}$$
From this condition, since $ f $ is odd, it follows that $ \frac{f(t)}{t^3} $ is decreasing in $ (-\infty , 0) $.

Moreover, in this case we can also weaken the condition $ (f_4) $, as follows.
$$ \mbox{ there  exists } q>4 \mbox{ \ and \ } M_1 > 0 \ \mbox{ such that}\  F(t) \geq M_1 |t|^q \ , \ \forall \ t\in \R . \leqno{(f_4 ')}$$

Our strategy to prove the second main result consists in applying the genus theory and it is inspired by the construction presented in \cite{[6]}.

\bt\label{t2}
Suppose $ (f_1 '), (f_2), (f_3), (f_4 '), (f_5) $. Then, problem (\ref{P}) admits a sequence of solution pairs $ \pm \un \in X $ such that $ I(\un) \RA + \infty $ as $ n \RA + \infty $. 
\et

The reader should pay attention that many difficulties arise while dealing with Choquard logarithmic equations with exponential nonlinearities. Concerning the logarithmic term, one will see that it is very difficult to guarantee when a Cerami sequence has a subsequence that converges strongly in $ X $ (defined in (\ref{eq33})). Inspired by \cite{[6]}, we obtain such convergence module translations. But, since the norm of $ X $ is not invariant under translations, new difficulties arise. Also, we can mention that the key convergence result that we have so far, does not work for functions satisfying $ \un(x) \RA 0 $ a.e. in $ \R $. Besides that, we have the difficulties concerning the controlling of exponential term, as we will mention later.

Next we make a quick overview of literature. In the first half of this overview, we recall that problems with nonlocal operators arise in many areas, such as optimization, finance, phase transitions, stratified materials, anomalous diffusion, crystal dislocation, soft thin films, semipermeable membranes, flame propagation, conservation laws and water waves. See e.g. \cite{[hitch],[caffarelli]}. For fractional problems involving $ (-\Delta)^s $, with $ N > 2s $ and $ s\in (0, 1) $, without logarithmic kernel, we refer to \cite{[r1], [r3]}, where the authors have obtained ground state solution under subcritical polynomial growth and, in addition, in \cite{[r3]}, they studied the regularity and derived some properties for those solutions. Moreover, in \cite{[r4]}, authors dealt with coercive potentials, while in \cite{[r2]}, they study the equation with potentials that vanish at infinity. In \cite{[frac]}, the authors guarantee the existence of ground state solution when the nonlinearity $ f $ has maximal exponential growth. We also refer the reader for the works \cite{[9], [weyl], [pucci1], [pucci2]}, in which the authors dealt with fractional Laplacian operator, and the works \cite{[5], [Lam], [17]}, for general problems with Moser-Trudinger type behaviour.

For the second half, we take a look into works that deal with Choquard logarithmic equations. In this sense, we can cite the recent works of \cite{[6], [cjj], [10], [wen], [alves]}, where the authors study problem (\ref{P}) in the local situation. 

In \cite{[6]}, the authors have proved the existence of infinitely many geometrically distinct solutions and a ground state solution, considering $ V: \Rdois \RA (0, \infty) $ continuous and $ \mathbb{Z}^2 $-periodic, $ \gamma > 0 $, $ b\geq 0 $, $ p\geq 4 $ and $ f(u)=b|u|^{p-2}u $. Here  because of the  periodic setting, the global Palais–Smale condition can fail, since  the corresponding functional become invariant under $ \mathbb{Z}^2 $-translations. Then, intending to fill the gap, in \cite{[10]} it is proved the existence of a mountain pass solution and a ground state solution for the local problem (\ref{P}) in the case $ V(x)\equiv a > 0 $, $ \lambda >0 $, $ 2<p<4 $ and $ f(u)=|u|^{p-2}u $. Also, they verified that, if $ p\geq 3 $, both levels are equal and provided a characterization for them. Moreover, in \cite{[cjj]}, the authors dealt with the existence of stationary waves with prescribed norm considering $ \lambda \in \R $ in a similar setting. Then, following the ideas of the predecessors, \cite{[wen]}, worked with local problem (\ref{P}) in $ \mathbb{R}^3 $ with a nonlinearity with polynomial growth. Finally, in \cite{[alves]}, the authors proved the existence of a ground state solution for local problem (\ref{P}), with a nonlinearity of Moser-Trudinger type. We also refer to \cite{[guo], [4], [14]} for more results about Choquard equations.

The present work aims to extend or complement those results already found in the literature, combining the fractional Laplacian operator with Choquard logarithmic equations and exponential nonlinearities. 

Throughout the paper, we will use the following notations: $ \Lt $ denotes the usual Lebesgue space with norm $ ||\cdot ||_t $ \ ; \ $ X' $ denotes the dual space of $ X $ \ ; \ $ B_r(x)=(x-r , x+r) $ is the ball in $ \R $ centred in $ x $ with radius $ r>0 $ and simply $ B_r $ when $ x=0 $ \ ; \ $ A^{c} $ stands for $ \R \setminus A $, for any subset $ A\subset \R $ \ ; \ $ C, C_1, C_2, ... $ will denote different positive constants whose exact values are not essential to the exposition of arguments. 

The paper is organized as follows: in section 2 we present the framework's problem and some technical and essential results, some of them already derived in previous works and whose application to our problem is immediate. Section 3 consists in the proof of a key proposition and our first main result. Finally, in section 4, we prove our second main result, considering the subcritical case. 

\section{Preliminary Results}

In this section we present the reader for the framework necessary to study problem (\ref{P}) and provide some technical results. 

We start remembering that the operator $ (-\Delta)^{\frac{1}{2}}: \mathcal{S}(\R) \RA \Ldois $ is given by
$$
(-\Delta)^{\frac{1}{2}}u(x) = C\left(1, \frac{1}{2}\right) \lim\limits_{\varepsilon \RA 0} \ds_{\R \setminus B_\varepsilon (x)}\dfrac{u(x)-u(y)}{|x-y|^2} dy \ , \ \forall \ x \in \R ,
$$
where the normalizing constant $ C\left(1, \frac{1}{2}\right) $ is defined in \cite{[hitch]} and $\mathcal{S}(\R)  $ denotes the Schwartz space. Equivalently, for $ u\in \mathcal{S}(\R) $, by \cite[Proposition 3.3]{[hitch]}, if $ \mathcal{F} $ denotes the Fourier Transform, 
$$
(-\Delta)^{\frac{1}{2}}u = \mathcal{F}^{-1}(|\xi|(\mathcal{F}u)) \ , \ \forall \ x \in \R .
$$
Moreover, in light of \cite[Proposition 3.6]{[hitch]}, we have
$$
||(-\Delta)^{\frac{1}{4}}u||_{2}^{2} = \dfrac{1}{2\pi} \intR \intR \dfrac{(u(x)-u(y))^2 }{|x-y|^2} dx dy \ , \ \forall \ u \in \Humeio ,
$$
even though we usually consider this equality omitting the normalizing constant $ \frac{1}{2\pi} $. 
 
Now, we turn our attention to the Hilbert space $$ W^{\frac{1}{2}, 2}(\R)=\Humeio = \left\{ u \in \Ldois \ ; \ \intR \intR \dfrac{|u(x)-u(y)|^2}{|x-y|^2} dx dy < + \infty \right\} , $$
endowed with the norm $ ||\cdot||^{2}=[\cdot]_{\frac{1}{2}, 2}^{2}+||\cdot||_{2}^{2} $, where
$$
[u]_{\frac{1}{2}, 2}^{2} = \intR \intR \dfrac{|u(x)-u(y)|^2}{|x-y|^2} dx dy .
$$
For further considerations about $ \Humeio $ and $ (-\Delta)^{\frac{1}{2}} $ and some useful results, we refer to \cite{[hitch], [demengel], [adams], [weyl], [pucci1], [pucci2]}.

Next, inspired in \cite{[21]}, in order to guarantee that the associated functional with the problem is well-defined, we consider the slightly smaller space
\begin{equation}\label{eq33}
X = \left\{ u\in \Humeio \ ; \ \intR \ln (1+|x|) u^2 (x) dx < + \infty \right\} .
\end{equation}
The space $ X $ endowed with the norm $ ||\cdot||_{X}^{2}=||\cdot||^{2}+||\cdot||_{\ast}^{2} $, where
$$
||u||_{\ast}^{2} = \intR \ln (1+|x|) u^2 (x) dx ,
$$
is a Hilbert space. 

Inspired by \cite{[6]}, we define three auxiliar symmetric bilinear forms 
\beqa
(u , v) \mapsto B_1(u, v)=\intR \intR \ln(1+|x-y|)u(x)v(y) dx dy ,
\eeqa
\beqa
(u , v) \mapsto B_2(u, v)=\intR \intR \ln\left(1+\dfrac{1}{|x-y|}\right)u(x)v(y) dx dy ,
\eeqa
\beqa
(u , v) \mapsto B_0(u, v)=B_1(u, v)-B_2(u, v)=\intR \intR \ln(|x-y|)u(x)v(y) dx dy. 
\eeqa
These definitions are understood to being over measurable function $u, v: \R \RA \R $, such that the integrals are defined in the Lebesgue sense. We also define the functionals  $V_1: \Humeio \RA [0, \infty],$ $V_2: L^{4}(\R) \RA [0, \infty) $ and $V_0: \Humeio \RA \R \cup \{\infty\},$ given by 
$V_1(u)=B_1(u^2 , u^2), $ $V_2(u)=B_2(u^2 , u^2) $ and $V_0(u)=B_0(u^2 , u^2),$ respectively.

By Hardy-Littlewood-Sobolev inequality (HLS) \cite{[15]}, using that $ 0 \leq \ln(1+r)\leq r $, for $ r>0 $, one can easily see that
\begin{align} \label{b2}
|B_2(u, v)| 
& \leq \intR \intR \dfrac{1}{|x-y|}u(x)v(y) dx dy 
 \leq K_0 ||u||_{2}||v||_{2} \ , \ \forall \  u, v \in L^{2}(\R) ,
\end{align}
where $ K_{0}>0 $ is the HLS constant.

As a consequence of (\ref{b2}), 
\beq \label{v2}
|V_2(u)|\leq K_0||u||_{4}^{4} \ , \ \ \forall \ u\in L^{4}(\R) ,
\eeq
so $ V_2 $ takes finite values over $ L^{4}(\R)$. Also, observing that
\beq \label{eq3}
\ln(1+|x\pm y|) \leq \ln(1+|x|+|y|) \leq \ln (1+|x|)+ \ln(1+|y|), \textrm{ \ for \ } x, y\in \R,
\eeq
we can estimate, applying Hölder Inequality,
\beq\label{b1}
B_1(uv, wz) \leq ||u||_{\ast}||v||_{\ast}||w||_2 ||z||_2 + ||u||_2 ||v||_2 ||w||_{\ast}||z||_{\ast} ,
\eeq
for all $ u, v, w, z \in \Ldois $. 

In our first lemma we determine in which spaces X can be embedded either continuously or compactly.

\bl\label{l4}
$X$ is continuously embedded in $ \Humeio $ and compactly embedded in $ \Lt $, for all $ t\in [2, +\infty) $. 
\el
\begin{proof}
The continuous embedding is immediate, since $ ||u||\leq ||u||_X $, for all $ u\in X $. In order to prove the compact embeddings, let $ (\un)\subset X $ such that $ \un \CF 0 $ in $ X $. Then, the result follows from the coerciveness of $\ln(1+|x|)$, \cite[Theorem 7.41]{[adams]}, \cite[Theorem 7.1]{[hitch]}, the interpolation inequality for $\Lt$ spaces and a diagonal argument.
\end{proof}

Next we make some considerations concerning the exponential behaviour of the nonlinearity. First of all, we mention the celebrated Moser-Trudinger Lemma \cite{[5], [ozawa], [frac]} 

\bl[Moser-Trudinger, \cite{[ozawa]}] \label{l5}
There exists $ 0<\omega \leq \pi $ such that, for all $ \alpha \in (0, \omega) $, there exists a constant $ H_\alpha >0 $ satisfying
$$
\intR (e^{\alpha u^2} -1) dx \leq H_\alpha ||u||_{2}^{2} ,
$$
for all $ u\in \Humeio $ with $ ||(-\Delta)^{\frac{1}{4}}u||_{2}^{2}\leq 1 $. 
\el
Based on the above lemma, \cite{[frac]} proved the following result.
\bl\label{l6} 
\cite[Proposition 2.1]{[frac]} For any $ \alpha >0 $ and $ u\in \Humeio $,
$$
\intR (e^{\alpha u^2}-1)dx < + \infty .
$$
\el
Hence, from Lemma \ref{l6} and equation (\ref{eq26}), for any $ u\in X $ we have
\begin{equation}\label{eq30}
\intR F(u) dx \leq \dfrac{1}{2}||u||_{2}^{2}+b_2||u||_{qr_2}^{q}\left(\intR (e^{r_1 \alpha u^2}-1) dx \right)^{\frac{1}{r_1}} < + \infty .
\end{equation}

The following lemma plays a key role on the continuity of the associated functional and while guaranteeing that the functional is lower semicontinuous for $ \Humeio $. In order to do that, we need a result that allow us to control the exponential term. 

\bl\label{l23}
(\cite[Proposition 2.7]{[n1]}) Let $ (\un) \subset X $ such that $ (\un) $ is strongly convergent in $ \Humeio $. Then, there exists a subsequence $ (\unk) \subset (\un) $ and a function $ g\in \Humeio $ satisfying $ |\unk(x)|\leq g(x) $ a.e. in $ \R $, for all $ k \in \N $. 
\el
\begin{proof}
The proof can be done similarly as in \cite{[n1]}. We only highlight here why the function $ w $ obtained in \cite{[n1]} belongs to $ \Humeio $. Since $ (w_n) $ is bounded in $ \Humeio $, passing to a subsequence, if necessary, $ w_n \CF v \in \Humeio $. Then, by \cite[Theorem 7.1]{[hitch]}, $ w_n \RA v $ in $ L^2(B_R) $, for all $ R>0 $. Hence, up to a subsequence, $ w_n(x)\RA v(x) $ a.e. in $ \R $. Once $ w_n(x)\RA w(x) $ a.e. in $ \R $, we conclude that $ w = v  $ a.e. in $ \R $. Therefore, $ w\in \Humeio $.
\end{proof}

\bl\label{l22}
Let $ (\un) \subset X $ and $ u\in X $ such that $ \un \RA u $ on $ \Humeio $. Then, we have
$$
\intR F(\un) \RA \intR F(u) \ \ , \ \ \intR f(\un)\un \RA \intR f(u)u \mbox{ \ \ and \ \ } \intR f(\un) v \RA \intR f(u)v \ , \ \forall \ v \in X .
$$
\el
\begin{proof}
Since $ \un \RA u $  in $ \Humeio $, $ \un(x) \RA u(x) $ a.e. in $ \R $ and, from \cite[Theorem 6.9]{[hitch]}, $ \un \RA u $ in $ \Lt $ for all $ t \geq 2 $. By Lemma \ref{l23} and using the Dominated Convergence Theorem, the result follows.
\end{proof}

We need some technical lemmas in order to obtain estimatives and convergence when dealing with logarithmic parts. Once the proofs found in \cite{[6]} remains essentially the same, by only exchanging $ \Rdois $ for $ \R $, we omit it here to make the paper concise.

\bl \label{l1}
(\cite[Lemma 2.1]{[6]}) Let $(\un)$ be a sequence in $\Ldois $ and $ u\in \Ldois \setminus \{0\} $ such that $\un \RA u $ pointwise a.e. on $\R $. Moreover, let $(\vn)$ be a bounded sequence in $\Ldois$ such that
$$
\sup\limits_{n\in \N}B_1(\un^2, \vn^2) < \infty .
$$
Then, there exist $n_0\in \N$ and $ C>0 $ such that $ ||\un||_{\ast}< C $, for $ n\geq n_0 $. If, moreover, 
$$
B_1(\un^2, \vn^2)\RA 0 \textrm{ \ \ and \ \ } ||\vn||_{2}\RA 0, \textrm{ \ as \ } n\RA \infty ,
$$
then 
$$
||\vn||_{\ast}\RA 0 \textrm{ \ ,  as \ } n \RA \infty .
$$
\el

\bl\label{l2}
(\cite[Lemma 2.6]{[6]}) Let $ (\un) $, $ (\vn) $ and $ (\wn) $ be bounded sequences in $ X $ such that $ \un \rightharpoonup u $ in $ X $. Then, for every $ z\in X $, we have $ B_1 (\vn \wn \ , \ z(\un - u))\RA 0 $, as $ n\RA + \infty $.
\el

\bl \label{l3} (\cite[Lemma 2.2]{[6]})
(i) The functionals $ V_0, V_1, V_2 $ are of class $ C^1 $ on $ X $. Moreover, $ V_{i}'(u)(v)=4B_i(u^2, uv) $, for $ u, v \in X $ and $ i=0, 1, 2 $. \\
(ii) $ V_2 $ is continuous (in fact continuously differentiable) on $\Lquatro$ .\\
(iii) $ V_1 $ is weakly lower semicontinuous on $ \Humeio $. 
(iv) $ I $ is lower semicontinuous on $ \Humeio $.
\el

Therefore, from Lemmas \ref{l22} and \ref{l3}, (\ref{eq30}) and \cite[Lemma 2]{[pucci2]}, we conclude that $ I: X \RA \R $ given by
$$
I(u) = \dfrac{1}{2}||u||^2 + V_0(u) - \intR F(u) dx
$$
is well-defined in $ X $ and $ I\in C^1(X, \R) $. 

Next, let us verify that $ I $ has the mountain pass geometry, in order to obtain a Cerami sequence for the mountain pass level $ c_{mp} $. 

\bl\label{l7}
There exists $ \rho >0 $ such that
\begin{equation}\label{eq4}
m_\beta = \inf\{ I(u) \ ; \ u\in X \ , \ ||u||=\beta \} > 0 \ , \ \forall \ \beta \in (0, \rho] 
\end{equation}
and
\begin{equation}\label{eq5}
n_\beta = \inf\{ I'(u)(u) \ ; \ u\in X \ , \ ||u||=\beta \} > 0 \ , \ \forall \ \beta \in (0, \rho] .
\end{equation}
\el
\begin{proof}
First of all, we choose $ r_1, r_2 >1 $ such that $ r_1\sim 1 $ e $ r_2 > 2 $. Then, consider $ u\in X\setminus\{0\} $ such that $ r_1 \alpha||u||^2 < \omega $, in order to apply the exponential estimatives. Thus, from (\ref{eq6}), HLS, Lemma \ref{l5} and Sobolev embeddings, we have
$$
I(u)=\dfrac{1}{2}||u||^2+\dfrac{1}{4}V_0 (u) - \intR F(u) dx \geq \dfrac{||u||^2}{2}[1-\varepsilon - C_2 ||u||^2 - C_3 ||u||^{q - 2}] .
$$
Therefore, since $ q - 2 >0 $, we can choose $ \rho >0 $ sufficiently small, such that (\ref{eq4}) is valid. 

Similarly, as
$$
I'(u)(u)= ||u||^{2}+V_0(u)-\intR f(u)u dx \geq ||u||^2 [ 1 - \varepsilon - C_4||u||^2 -C_5 ||u||^{q -2}] ,
$$
once again one can pick $ \rho >0 $ sufficiently small, such that (\ref{eq5}) holds. 
\end{proof}

\bl\label{l8}
Let $ u\in X\setminus \{0\} $ and $ q>4 $. Then, 
$$
\lim\limits_{t\RA 0} I(tu) = 0 \ \ , \ \ \sup\limits_{t>0} I(tu) < + \infty \ \ \mbox{ and \ \ } I(tu)\RA - \infty \ , \ \mbox{ as \ } t\RA + \infty .
$$
\el
\begin{proof}
Let $ u\in X\setminus \{0\} $. First, from $ (f_4) $ and $ q>4 $,
$$
I(tu) = \dfrac{t^2}{2}||u||^2 + \dfrac{t^{4}}{4}V_0(u) - \intR F(tu) \leq \dfrac{t^2}{2}||u||^2 + \dfrac{t^{4}}{4}V_0(u) - C_q t^q ||u||_{q}^{q} \RA - \infty ,
$$
as $ t\RA + \infty $. Now, from (\ref{eq26}) and Lemma \ref{l5}, taking $ t>0 $ sufficiently small such that $ r_1 \alpha t ||u||^2 < \omega $, we have
$$
\intR F(tu) dx \leq \dfrac{t^2}{2}||u||_{2}^{2} + H_\alpha t^{q}||u||_{qr_2}^{q} \RA 0 ,
$$
as $ t\RA 0 $. Hence, $ I(tu)\RA 0 $ as $ t\RA 0 $. Finally, this behaviour combined with the fact that $ I $ is $ C^1 $, tell us that $ \sup\limits_{t>0} I(tu) < + \infty $.
\end{proof}

For the next results, consider a sequence $ (\un) \subset X$ satisfying
\begin{equation}\label{eq7}
\exists \ d > 0 \mbox{ \ s.t. \ } I(\un) < d \ , \ \forall \ n\in \N \ ,  \mbox{ \ and \ } ||I'(\un)||_{X'}(1+||\un||_X) \RA 0 \mbox{ \ , as \ } n \RA + \infty
\end{equation}

\bl\label{l9}
Let $ (\un)\subset X $ satisfying (\ref{eq7}). Then, $ (\un) $ is bounded in $ \Humeio $.
\el
\begin{proof}
From (\ref{eq7}) and $ (f_3) $, we have
$$
d+ o(1) \geq I(\un) - \dfrac{1}{4} I'(\un)(\un) \geq \dfrac{1}{4}||\un||^2 + \left(\dfrac{\theta}{4}-1\right) \intR F(\un) dx \geq \dfrac{1}{4}||\un||^2 .
$$
Hence, we conclude that $ (\un) $ is bounded in $ \Humeio $.
\end{proof}

\bo\label{obs1}
One can easily verify, from Lemma \ref{l8} and the Intermediate Value Theorem, that $ 0 < m_\rho \leq c_{mp} < +\infty $. Then, since $ I $ has the mountain pass geometry, similarly as in \cite[Lemma 3.2]{[10]}, there exists a sequence $ (\un)\subset X $ such that
\begin{equation}\label{eq8}
I(\un) \RA c_{mp} \mbox{ \ \ \ and \ \ \ } ||I'(\un)||_{X'}(1+||\un||_X) \RA 0 \mbox{ \ , as \ } n \RA + \infty .
\end{equation}
Moreover, this sequence satisfies (\ref{eq7}).
\eo

Subsequently, we will guarantee that the norms of the sequence obtained in Remark \ref{obs1} can be taken sufficiently small such that the exponential estimatives are valid for it. Moreover, this will be valid for any sequence $ (\un) \subset X $ satisfying $ I(\un) \leq c_{mp} $, for all $ n\in \N $.

\bl\label{l10}
Let $ (\un)\subset X $ satisfying (\ref{eq8}) and $ q>4 $. Then, for some $ \rho_0 >0 $ sufficiently small,
$$
\limsup\limits_{n} ||\un||^2 < \rho_{0}^{2} .
$$
\el
\begin{proof}
From Lemma \ref{l9}, we have that $ (\un) $ is bounded in $ \Humeio $ and $ 4c_{mp} + o(1)\geq ||\un||^2 $, for all $ n\in \N $. Then, $ \limsup\limits_{n}||\un||^2 \leq 4c_{mp} $. Hence, we need to find an estimative for the value $ c_{mp} $. To do that, consider the set $ \A = \{u \in X \ ; \ u\neq 0 \ , \ V_0(u) \leq 0 \} $. Defining $ u_t(x)=t^2u(tx) $, for all $ t>0 $, $ u\in X\setminus\{0\} $ and $ x\in \R $, we verify that
$$
V_0(u_t)=t^6 V_0 (u) - t^6 \ln t ||u||_{2}^{4} \RA - \infty ,
$$
as $ t\RA + \infty $, and hence $ \A \neq \emptyset $. Moreover, from immersions \cite[Theorem 6.9]{[hitch]}, there exists $ C>0 $ such that $ ||u||\geq C||u||_q $, for all $ u\in \Humeio\setminus\{0\} $. Thus, it makes sense to define
$$
S_q(v)=\dfrac{||v||}{||v||_q} \mbox{ \ \ \ and \ \ \ } S_q = \inf\limits_{v\in \A} S_q(v) \geq \inf\limits_{v\neq 0}S_q(v) > 0 .
$$
Now we are ready to estimate $ c_{mp} $. From Lemma \ref{l8}, for $ v\in \A $ and $ T>0 $ sufficiently large, $ I(Tv)<0 $. So, we can define $ \gamma\in \Gamma $ by $ \gamma(t)=tTv $, such that
$$
c_{mp} \leq \max\limits_{0 \leq t \leq 1}I(\gamma(t)) = \max\limits_{0 \leq t \leq 1}I(tTv) \leq \max\limits_{t\geq 0}I(tv) .
$$
Consequently, for $ \psi \in \A $, we have
$$
c_{mp}\leq \max\limits_{t\geq 0}I(t\psi) \leq \max\limits_{t\geq 0}\left\{ \dfrac{t^2}{2}||\psi||^2-C_qt^q||\psi||_{q}^{q}\right\} \leq \left(\dfrac{1}{2}-\dfrac{1}{q}\right)\dfrac{S_q(\psi)^{\frac{2q}{q-2}}}{(qC_q)^{\frac{2}{q-2}}} .
$$
Thus, taking the infimum over $ \psi \in \A $, we obtain
$$
\limsup\limits_{n}||\un||^2 \leq \dfrac{2(q-2)}{q}\dfrac{S_{q}^{\frac{2q}{q-2}}}{(qC_q)^{\frac{2}{q-2}}} \leq \rho_{0}^{2} ,
$$
for $ C_q >0 $ sufficiently large.
\end{proof}

\bl\label{l11}
Let $ (\un)\subset X $ be bounded in $ \Humeio $ such that 
\begin{equation}\label{eq9}
\liminf \sup\limits_{y\in \mathbb{Z}} \ds_{B_2(y)}|\un(y)|^2 dx > 0 .
\end{equation}
Then, there exists $ u\in \Humeio \setminus\{0\} $ and $ (\yn)\subset \mathbb{Z} $ such that, up to a subsequence, $ \yn \ast \un = \until \CF u\in \Humeio $.  
\el
\begin{proof}
From (\ref{eq9}), $ \liminf $ properties, the boundedness of $ (\un) $ in $ \Humeio $, \cite[Theorem 7.41]{[adams]} and  \cite[Theorem 7.1]{[hitch]}, one can construct the desired subsequence and obtain the result.
\end{proof}

\section{Proof of Theorem \ref{t1}}

In the present section, we finish the proof of Theorem \ref{t1}. As a first step, we prove a key proposition which provide us with nontrivial critical points.

\bp\label{p1}
Let $ q>4 $ and $ (\un)\subset X $ satisfying (\ref{eq8}). Then, passing to a subsequence if necessary, only one between the following items is true:
\begin{itemize}
\item[(a)] $ ||\un||\RA 0 $ and $ I(\un) \RA 0 $, as $ n \RA + \infty $.

\item[(b)] There exists points $ \yn \in \mathbb{Z} $ such that $ \until = \yn \ast \un \RA u $ in $ X $, for a nontrivial critical point $ u\in X $ for $ I $.
\end{itemize}
\ep
\begin{proof}
First of all, from Lemma \ref{l9}, $ (\un) $ is bounded in $ \Humeio $ and, passing to a subsequence if necessary, by Lemma \ref{l10}, 
$$
\intR (e^{r_1 \alpha \un^{2}}-1 ) dx \leq H_\alpha \ , \ \forall \ n \in \N .
$$
Suppose that (a) does not happen. 

\noindent \textbf{Claim 1:} $ \liminf \sup\limits_{y\in \mathbb{Z}} \ds_{B_2(y)}|\un(x)|^2 dx > 0  $. 

Suppose the contrary. Then, from an easy adaptation of Lion's Lemma found in \cite[Lemma 2.4]{[yu]}, one has $ \un \RA 0 $ in $ L^{t}(\R) $, for all $ t>2 $. Thus, from (\ref{v2}), $ V_2(\un) \RA 0 $. Moreover, since $ q > 2 $, from (\ref{eq7}),
$$
\left|\intR f(\un) \un dx \right| \leq \varepsilon ||\un||_{2}^{2} + C_1 ||\un||_{qr_2}^{q} \RA 0 ,
$$
as $ \varepsilon \RA 0 $ and $ n \RA + \infty $. Consequently,
$$
||\un||^2 +V_1(\un) = I'(\un)(\un) +V_2(\un) + \intR f(\un)\un dx \RA 0 ,
$$
as $ n\RA + \infty $. Thus, $ ||\un||\RA 0 $ and $ V_1(\un) \RA 0 $. From the continuous immersions \cite[Theorem 6.9]{[hitch]}, $ ||\un||_2 \RA 0 $. Then, from (\ref{eq6}), $ \intR F(\un) dx \RA 0 $. As a consequence, $ I(\un) \RA 0 $, which is a contradiction. Therefore, the claim is valid. 

From Claim 1 and Lemma \ref{l11}, up to a subsequence, there exists $ (\yn)\subset \mathbb{Z} $ and $ u\in \Humeio \setminus\{0\} $ such that $ \until= \yn \ast \un \CF u $ in $ \Humeio $. We can conclude that $ (\until) $ is bounded in $ \Lt $, for all $ t\geq 2 $, and $ \until(x)\RA u(x) $ a.e. in $ \R $. 

Observe that, as $ q > 2 $, 
\begin{align*}
V_1(\until)=V_1(\un)& =I'(\un)(\un) +V_2(\un) + \intR f(\un)\un dx -||\un||^2\\
& \leq K_0 ||\un||_{4}^{2}+\varepsilon C_2 + C_1 ||\un||_{qr_2}^{q} + o(1) \leq C_3 + o(1) .
\end{align*}
That is, $ \sup\limits_{n}V_1(\until) < + \infty $. So, from $ (\until) $ being bounded in $ \Ldois $ and Lemma \ref{l1}, $ (||\until||_\ast) $ is bounded. Therefore, $ (\until) $ is bounded in $ X $ and, once $ X $ is reflexive, passing to a subsequence, if necessary, $ \until \CF u $ in $ X $. From Lemma \ref{l4}, $ \until \RA u $ in $ \Lt $, for all $ t\geq 2 $. 

\noindent \textbf{Claim 2:} $ I'(\until)(\until - u)\RA 0 $, as $ n\RA + \infty $.

One can easily see, by a change of variable, that $ I'(\until)(\until - u) = I'(\un)(\un - (-\yn)\ast u) $. Thus,
\begin{equation}\label{eq48}
| I'(\until)(\until - u)| = |I'(\un)(\un - (-\yn)\ast u)| \leq ||I'(\un)||_{X'}(||\un||_X + ||(-\yn)\ast u||_X) \ , \ \forall \ n\in \N .
\end{equation}
Then, we first seek for an useful inequality for $ ||(-\yn) \ast u||_X $. Once $ ||\un||_X $ already appears in (\ref{eq48}), a natural way is to look for a constant $ C>0 $ such that $||(-\yn)\ast u||_X \leq C ||\un||_X  $. 

We start analysing $ ||\un||_\ast $. If $|\yn|\RA +\infty$, then, for $ x\in \R $,
$$
\ln(1+|x-\yn|)-\ln(1+|\yn|)=\ln \left(\dfrac{1+|x-\yn|}{1+|\yn|}\right) \RA 0 , n\RA + \infty .
$$
Therefore, there exists $ C_4 > 0 $ such that $ \ln(1+|x-\yn|)\geq C_4 \ln(1+|\yn|) $. \\
Now, suppose that $ (\yn)\subset \mathbb{Z} $ converges to $ y_0 \in \mathbb{Z} $. Then, up to a subsequence, $ \yn \equiv y_0 $. Suppose $ y_0 > 0 $. Let $ \delta > 0 $ such that $ \delta < |y_0| $ and define $ \Omega = (-\delta , 0) $. Hence, for $ x\in \Omega $, we have that $ |x-y_0| > |y_0| $. Therefore, by the Mean Value Theorem, there exists $ x_\delta \in \Omega $, satisfying
\begin{align*}
||\un||_{\ast}^{2} & \geq \int\limits_{\Omega} \ln(1+|x-\yn|)\until^2 (x) dx \\
& = \delta \ln(1+|x_\delta -\yn|)\until^2 (x_\delta ) \\
& = C_4 \ln(1+|x_\delta -y_0|) \geq C_4 \ln(1+|y_0|) = C_4 \ln(1+|\yn|) ,
\end{align*}
for $ C_4 > 0 $. Analogously for $ y_0 < 0 $.  If $ y_0 =0 $, $ \until = \un $ and nothing remains to be proved. So, in any of the cases, there exists $ C_4 >0 $ such that
\begin{equation}\label{eq47}
||\un||_{\ast}^{2}= \intR \ln(1+|x -\yn|)\until^2 (x) dx \geq C_4 \ln(1+ |\yn|) \ , \forall \ n \in \N .
\end{equation}
Now, from (\ref{eq3}), we have
$$
||\until||_{\ast}^{2} = \intR \ln (1+|x+\yn|) \un^2(x) dx \leq ||\un||_{\ast}^{2} + \ln (1+|\yn|) ||\un||_{2}^{2} .
$$
From this, (\ref{eq3}) and (\ref{eq47}), since every norm is weakly lower semicontinuos, $ \until \rightharpoonup u $ in $ X $ and $ ||\cdot||_2 $ is $ \mathbb{Z}$-invariant, follows that
\begin{align*}
||(-\yn) \ast u||_{\ast}^{2}& \leq ||u||_{\ast}^{2}+\ln(1+|\yn|)||u||_{2}^{2} \\
& \leq ||\until||_{X}^{2} + \ln(1+|\yn|)||\un||_{2}^{2} \\
& = ||\un||^{2}  + ||\un||_{\ast}^{2} (1 + C_5||\un||_{2}^{2}) \\
& \leq ||\un||^{2}  + C_6 ||\un||_{\ast}^{2} \leq C_7 ||\un||_{X}^{2}
\end{align*}
for $ n\in \N $ and $ C_7 >0 $. Consequently, there exists a constant $ C_8 >0 $ such that, after passing to a subsequence, we have, for all $ n\in \N $,
$$
||(-\yn) \ast u||_X^{2} =||u||^2 + ||(-\yn) \ast u||_{\ast}^{2} \leq ||\un||^2 +C_7 ||\un||_{X}^{2} \leq C_8 ||\un||_{X}^{2} .
$$
Thus, from (\ref{eq8}),
$$
| I'(\until)(\until - u)|\leq (1+C_{8}^{\frac{1}{2}})||I'(\un)||_{X'}||\un||_X \RA 0 \ , \ \mbox{as \ } n \RA + \infty .
$$

\noindent \textbf{Claim 3:} $ \intR f(\until)(\until - u ) dx \RA 0 $, as $ n\RA + \infty $.

Since $ ||\cdot|| $ is $ \mathbb{Z} $-invariant, Moser-Trudinger inequality and the exponential results above hold for the sequence $ (\until) $ as well. It follows that
$$
\left| \intR f(\until)(\until -u) dx \right| \leq ||\until||_{2}||\until -u||_2 + C_9 ||\until -u||_{qr_2}^{q} \RA 0 ,
$$ 
as $ n \RA + \infty $, concluding the claim.

Moreover, observe that $$
\left|V_{2}'(\until)(\until - u)\right| \leq K_0 ||\until||_{4}^{3}||\until -u||_4  \RA 0
$$
and
$$
V_{1}'(\until)(\until - u)=B_1 (\until^2 , \until (\until - u)) = B_1 (\until^2 ,  (\until - u)^2) + B_1(\until^2 , u (\until - u)) .
$$
Also, since $ (\until) $ is bounded in $ X $, from Lemma \ref{l2}, $ B_1(\until^2 , u (\until - u))\RA 0 $. So, from the above estimations 
$$
o(1) =||\until||^{2}-||u||^2 +  B_1 (\until^2 ,  (\until - u)^2) + o(1) \geq ||\until||^{2}-||u||^2  + o(1) .
$$
Hence, $ ||\until||\RA ||u|| $ and $B_1 (\until^2 ,  (\until - u)^2)\RA 0$. Therefore, we conclude that $ ||\until - u||\RA 0 $ and, combined with Lemma \ref{l1}, we have $ ||\until - u ||_{\ast}\RA 0 $. Consequently, $ \until \RA u $ in $ X $. 

Finally, we have to show that $ u $ is a critical point of $ I $. Let $ v\in X $. So, as we deduced above, it is possible to find $ C_{10} >0 $ such that
$$
||(-\yn) \ast v||_X \leq C_{10} ||\un||_X \ , \ \forall \ n \in \N .
$$
Thus, 
$$
|I'(u)(v)| = \lim\limits_{n\RA + \infty} |I'(\until)(v)|=\lim\limits_{n\RA + \infty} |I'(\un)((-\yn)\ast v)| \leq C_{10} \lim\limits_{n\RA + \infty} ||I'(\un)||_{X'}||\un||_X = 0 .
$$
\end{proof}

\begin{proof}[Proof of Theorem \ref{t1}]
\textbf{(i)} From Lemma \ref{l7} and Proposition \ref{p1} there exists a nontrivial critical point of $ I $, $ u_0 \in X $, such that $ I(u_0)=c_{mp} $. 

\noindent \textbf{(ii)} Define $ K =\{v \in X\setminus \{0\} \ ; \ I'(v)=0\} $. Since $ u_0\in K $, $ K\neq \emptyset $. Thus, we can consider $ (\un)\subset K $ satisfying $ I(\un)\RA c_g = \inf\limits_{v\in K} I'(v) $.

Observe that $ c_g \in [-\infty , c_{mp}] $. If $ c_g = c_{mp} $, nothing remains to be proved. Otherwise, if $ c_g < c_{mp} $, passing to a subsequence if necessary, we can assume that $ I(\un) \leq c_{mp} $, for all $ n \in \N $. So, Lemma \ref{l10} holds for this sequence and we can apply all the results derived earlier. Therefore, from Lemma \ref{l7}, $ (\un)\subset K $ and Proposition \ref{p1} there exists $ (\yn)\subset \mathbb{Z} $ such that $ \until \RA u $ in $ X $, for a nontrivial critical point $ u $ of $ I $ in $ X $. Consequently, 
$
I'(u)=\lim I'(\until)=\lim I'(\un)=0 
$
and we conclude that $ u\in K $ and 
$
I(u)=\lim I(\until) = \lim I(\un) = c_g .
$
Particularly, we see that $ c_g > -\infty $.
\end{proof}

\section{Proof of Theorem \ref{t2}}

To finish this section, we will verify some properties of an important auxiliary function, namely, $ \varphi_u : \R \RA \R $, given by $ \varphi_u(t) = I(tu) $, for all $ u\in X\setminus\{0\} $ and $ t \in \R $.

\bl\label{l110}
Let $ u\in X\setminus \{0\} $. Then, $ \varphi_u $ is even and there exists a unique $ t_u \in (0, +\infty) $ such that $ \varphi_u '(t)>0 $, for all $ t\in (0, t_u) $, and $ \varphi_u '(t)<0 $, for all $ t\in (t_u , \infty) $. Moreover, $ \varphi_u (t) \RA -\infty $, as $ t \RA +\infty $.
\el
\begin{proof}
Since $ f $ is odd, $ I $ is even and, consequently, $ \varphi_u $ is even as well. 

\noindent \textbf{(i)} For $ t > 0 $ sufficiently small and $ \alpha > 0 $, for (\ref{eq27}), we have
$$
\varphi_u '(t) \geq t||u||^2[1- C_2t^2 ||u||^2 -\varepsilon -C_4t^{q-2}||u||^{q-2}] .
$$
Thus, $ \varphi_u '(t) > 0 $ for $ t> 0 $ sufficiently small. 

\noindent \textbf{(ii)} From $ (f_3) $ and $ (f_4 ') $, once $ q>4 $, 
$$
\varphi_u '(t)  \leq t||u||^2 + t^3 V_1(u) - C_3 t^{q-1}||u||_{q}^{q} \RA - \infty \mbox{ \ , as \ } t\RA + \infty .
$$

Hence, from (i)-(ii), since $ I $ is $ C^1 $, there exists $ t_u \in (0, + \infty) $ such that $ \varphi_u '(t_u) =0 $, which is unique by $ (f_5) $. 
\end{proof}

\bc\label{c12}
Let $ u\in X\setminus \{0\} $. Then, there exists a unique $ t_{u}'\in (0, + \infty )$ such that $ \varphi_u(t)>0 $, for $ t\in (0, t_{u}') $, and $ \varphi_u(t)<0 $, for $ t\in (t_{u}', + \infty) $. Moreover, $ t_u $ given by Lemma \ref{l110} is a global maximum for $ \varphi_u $.
\ec

\bl\label{l111}
For each $ u\in X\setminus \{0\} $, the map $ u \mapsto t_u ' $ is continuous.
\el
\begin{proof}
First of all, observe that, from the uniqueness of $ t_u ' $, the map is well-defined for all $ u\in X\setminus \{0\} $. Consider now $ (\un) \subset X \setminus \{0\} $ such that $ \un \RA u $ in $ X\setminus\{0\} $. 

Suppose, by contradiction, that $ (t_{\un}') $ is not bounded. So, from Lemma \ref{l110}, the definition of $ \varphi_u $, the continuity of $ I $ and Corollary \ref{c12}, there exists $ t_0 > 0 $ such that $ 0 \leq \varphi_u (t_0)<0 $, which is a contradiction. Therefore, $ (t_{\un}') $ is bounded. 

Hence, there exists $ (t_{\unk}')\subset (t_{\un}') $ such that $ t_{\unk} ' \RA \overline{t}\in (0, + \infty) $. Thus, we infer that
$$
0  \leq ||t_{\unk} ' \unk - \overline{t}u||_X  \leq |t_{\unk}' - \overline{t}| ||\unk||_X + \overline{t}||\unk - u||_X \RA 0 .
$$
Consequently, $ t_{\unk}'\unk \RA \overline{t}u $ in $ X $. From continuity of $ I $, $ I(t_{\unk}'\unk ) \RA I(\overline{t}u) $. 

On the other hand, from $ t_{\unk}' $ definition, $ I(t_{\unk}'\unk)=0 $. Hence, $ I(\overline{t}u)=0 $.  From uniqueness given by Corollary \ref{c12}, $ \overline{t}=t_u ' $. So, we conclude that $ t_{\unk}'\RA t_u ' $. Since this argument can be done to any subsequence of $ (t_{\un}) $, by standard arguments one can see that $ t_{\un}' \RA t_u ' $, proving  the continuity. 
\end{proof}

In this section we will provide the proof of the multiplicity result stated. In order to do so, we will need to verify some results concerned with the genus theory, denoted by $ \gamma $ and which definition, given over $ \A = \{ A\subset X \ ; \ A \mbox{ is symmetric and closed} \} $ (with respect to the continuity in $ X $), and basic properties can be found in \cite[Chapter II.5]{[genus]}. 

We start defining the sets 
$$
K_c = \{ u\in X \ ; \ I'(u)=0 , I(u)=c \} \ , \  c\in (0, + \infty)
$$
and
$$
A_{c, \rho} = \{ u\in X \ ; \ ||u-v|| \leq \rho \mbox{ \ , for some } v\in K_c\} \ , \  c\in (0, + \infty).
$$
It is easy to verify that the sets $ K_c $ and $ A_{c, \rho} $ are symmetric, closed (with respect to $ ||\cdot||_X $) and invariant under $ \mathbb{Z} $ translations, i.e, if $ u\in K_c , A_{c, \rho} $, then $ z\ast u \in K_c , A_{c, \rho} $, for all $ z\in \mathbb{Z} $.

Next, we fix a continuous map $ \beta: \Ldois \setminus \{0\} \RA \R $ that is equivariant under $ \mathbb{Z} $ translations, i.e, $ \beta(x \ast u) = x+\beta(u) $, for $ x\in \mathbb{Z} $ and $ u\in \Ldois \setminus \{0\} $.  We also require that $ \beta(-u)=\beta(u) $. Such map is called a \textit{generalized barycenter map} and an example can be constructed as that one in \cite{[bar]}. Hence, we can define the sets
$$
\tilde{K}_c = \{ u\in K_c \ ; \ \beta(u)\in [-4, 4] \} ,
$$
which are clearly symmetric. Moreover, before we enunciate our first result, we need to recall the Gauss bracket $ [ \cdot ] : \R \RA \mathbb{Z} $, given by $ [s]=\max\{n\in \mathbb{Z} \ ; \ n \leq s \} $ (see \cite[Chapter 3]{[n2]}. We recall some properties of the Gauss bracket, that are needed inside of the proof of our result.

\bl\label{l19}
Let $ [ \cdot ] : \R \RA \mathbb{Z} $, given by $ [s]=\max\{n\in \mathbb{Z} \ ; \ n \leq s \} $. Then,

\noindent \textbf{(i)} $ 0 \leq s - [s] < 1 $, for all $ s\in \R $.

\noindent \textbf{(ii)} if $ z\in \mathbb{Z} $ and $ s\in \R $, then $ [z+s]=z+[s] $.

\noindent \textbf{(iii)} Let $ s\in \R $ such that $ s - [s] \geq \dfrac{1}{2} $, then $ s - \dfrac{1}{2} - \left[ s - \dfrac{1}{2}\right] < \dfrac{1}{2} $.

\noindent \textbf{(iv)} $ 0 \leq s - [s] < 1 $, for all $ s\in \R $.

\noindent \textbf{(v)} Let $ s\in \R $. Then, $ s- \left[s - \dfrac{1}{2}\right]< 1 $.
\el

\bp\label{p19}
Let $ c>0 $. Then, there exists $ \rho_0 = \rho_0(c) > 0 $ such that $ \gamma(A_{c, \rho}) < \infty $, for all $ \rho \in (0, \rho_0) $.
\ep
\begin{proof}
Fix a baricenter map $ \beta : \Ldois \setminus \{0\} \RA \R $. 

\noindent \textbf{Claim 1:} The sets $ \tilde{K}_c $ are compact in $ X $.

Let $ (\un)\subset \tilde{K}_c $. From Proposition \ref{p1}, there exists $ (\yn)\subset \mathbb{Z} $ such that, passing to a subsequence if necessary, $ \wn = \yn \ast \un \RA u \in K_c $ in $ X $. Consequently, as $ \beta $ is continuous,
$$
\yn + \beta(\un) = \beta(\yn \ast \un) = \beta (\wn) \RA \beta(u) .
$$
Since $ (\un)\subset \tilde{K}_c $, $ \beta(\un)\subset [-4, 4] $. From compactness, passing to a subsequence if necessary, there exists $ r\in [-4, 4] $ such that $ \beta(\un)\RA r $ in $ \R $. Also, since $ (\beta(\un)) $ and $ (\yn + \beta(\un)) $ are bounded, $ (\yn)\subset \mathbb{Z} $ is bounded. Thus, we can assume that $ (\yn) $ is convergent. Therefore,
$$
\beta(u)=\lim \yn + \lim \beta(\un) \Rightarrow \lim \yn = \beta(u) - r .
$$
Moreover, since $ (\yn)\subset \mathbb{Z} $, $ \beta(u)-r=\overline{y} \in \mathbb{Z} $.

Lets prove that $ \un \RA u_0 = (-\overline{y})\ast u $, with respect to $ ||\cdot||_X $. First of all, observe that, since $ (\yn)\subset \mathbb{Z} $ and $ \yn \RA \overline{y} \in \mathbb{Z} $, there exists $ n_0\in \N $, such that for all $ n\geq n_0 $, $ \yn = \overline{y} $. Moreover, 
\begin{align*}
||\un - (-\yn)\ast u||_\ast & \leq \intR \ln(1+|y|)|\yn \ast \un - u|^2 dy + \intR \ln(1+|\yn|) |\yn \ast \un - u |^2 dy  \\
& \leq (1+ \ln(1+|\yn|))||\yn \ast \un -u||_{X}^{2} .
\end{align*}
Since $ \yn = \overline{y} $, for all $ n\geq n_0 $, and $ \yn \ast \un \RA u $ in $ X $, we have
\begin{align*}
\lim ||\un - (-\yn)\ast u||_{\ast}^{2} & \leq \lim [(1+ \ln(1+|\yn|))||\yn \ast \un -u||_{X}^{2} ] \\
& = (1+ \ln(1+|\overline{y}|))\lim |\yn \ast \un -u||_{X}^{2} = 0.
\end{align*}
Also, as $ ||\yn \ast \un - u||=||\un -(-\yn)\ast u|| $ and $ \yn \ast \un \RA u $ in $ \Humeio $, $ ||\un - (-\yn)\ast u|| \RA 0 $, as $ n \RA + \infty $. Hence, $ ||\un - (-\yn)\ast u||_X \RA 0 $. Finally, since $ \yn = \overline{y} $, for all $ n \geq n_0 $, we conclude that $ \un \RA u_0 $, as desired. 

Therefore, from the $ \mathbb{Z} $-invariance of $ K_c $, $ (-\overline{y})\ast u\in K_c $ and, by continuity of $ \beta $, $ \beta(\un)\RA \beta(u_0) \in [-4, 4] $. Consequently, $ u_0\in \tilde{K}_c $, and the claim is valid. 

For $ c>0 $, from the facts that $ \beta $ is even and $ X \hookrightarrow \Humeio $, one can see that $ \tilde{K}_c $ is compact in $ \Humeio $ and, also, $ \tilde{K}_c \in \A $. Thus, since $ 0 \not\in \tilde{K}_c $, $ \gamma(\tilde{K}_c)< + \infty $.

Hence, from the genus definition, there exists a continuous and odd function $ g_0 : \tilde{K}_c \RA \mathbb{R}^k \setminus \{0\} $, for some $ k \geq 0 $. Then, from Claim 1, one can easily verify that $ g_0 $ is continuous with respect to $ ||\cdot|| $.

Therefore, since $ \Humeio $ is a normal topological space and $ \tilde{K}_c $ is closed in $ \Humeio $, we can apply a corollary of Tietze's Theorem in order to extend $ g_0 $ to a continuous and odd function (with respect to $ ||\cdot|| $), $ g: \Humeio \RA \R $. 

\noindent \textbf{Claim 2:} There exists $ \overline{\rho} >0 $ such that $ g(u)\neq 0 $, for all $ u\in U_{\overline{\rho}}=\{u\in \Humeio \ ; \ ||u-v||\leq \overline{\rho} \mbox{ \ , for some \ } v\in \tilde{K}_c \}$.

First of all, observe that $ g(\tilde{K}_c)=g_0(\tilde{K}_c)\subset \mathbb{R}^k\setminus \{0\} $, that is, $ \{0\}\cap g(\tilde{K}_c) = \emptyset $ or, equivalently, $ g^{-1}(\{0\})\cap \tilde{K}_c = \emptyset $. From continuity of $ g $, we have that $g^{-1}(\{0\})  $ is closed in $ \Humeio $. Consequently, $d= dist(g^{-1}(\{0\}), \tilde{K}_c) >0 $. 

Set $ \overline{\rho}=\dfrac{d}{2} $. Suppose, by contradiction, that there exits $ u\in U_{\overline{\rho}} $ such that $ g(u)=0 $, i.e, $ u\in  g^{-1}(\{0\}) $. Thus, there exists $ v\in \tilde{K}_c $ satisfying $$ dist(u, v)=||u-v|| \leq \overline{\rho} = \dfrac{d}{2} .$$ 
Hence, $ d \leq ||u-v|| \leq \dfrac{d}{2} $, which is a contradiction. 

Therefore, $ g(u)\neq 0 $ for all $ u\in U_{\overline{\rho}} $.

Also, one can verify that $ U_{\overline{\rho}}\in \A $. In the following, we need to define the set
$$
L_1 = \left\{ u\in \Lp \setminus \{0\} \ ; \ |\beta(u)-[\beta(u)]|\leq \dfrac{1}{2}\right\} .
$$
One can easily see that $ z\ast u \in L_1 $, for all $ z\in \mathbb{Z} $ and $ u\in L_1 $. Now, define a map $ h_1 : L_1 \RA L_1 $ by $ h_1(u)=(-[\beta(u)])\ast u $. It is easy to verify that $ h_1 $ is well-defined, odd, invariant under $ \mathbb{Z} $ translations and an isometry. Also, setting $ a = \dfrac{1}{2}$ and $ L_2 = a \ast L_1 \subset \Ldois \setminus \{0\} $, it is easy to see, from item (iii) of Lemma \ref{l19}, that $ \Ldois \setminus \{0\} \subset L_1 \cup L_2  $.

Now, define the map $ h_2 : L_2 \RA L_2  $ by $ h_2(u)=a \ast [h_1 ((-a)\ast u)] $, for $ u\in L_2 $. Once again, it is easy to verify that $ h_2 $ is odd, $ \mathbb{Z} $-invariant and an isometry. Moreover, it follows immediately from the properties of $ \beta $, definition of $ h_2 $ and $ h_1 $ and Lemma \ref{l19}-(v), that $ \beta(h_2(u))\in [0, 1] $ for  $ u\in L_2 $.
 
Next, define the sets $ A_i = h_{i}^{-1}(U_{\overline{\rho}})\subset L_i $, for $ i= 1, 2 $. One can see that $ A_i \in \A $ and that they are $ \mathbb{Z} $-invariant, for $ i=1, 2$. Set $ A = A_1 \cup A_2 $.

\noindent \textbf{Claim 3:} There exists $ \rho_c \in (0, +\infty) $, such that $ A_{c, \rho}\subset A $, for all $ \rho \in (0, \rho_c) $.

Suppose, by contradiction, that for every $ n\in \N $, there exists $ \rho_n \in (0, \frac{1}{n}) $ such that $ A_{c, \rho_n}\nsubseteq A $. Then, there exists $ \un\in A_{c, \rho_n}\setminus A $, for every $ n\in \N $. 

We can suppose that $ \un \in L_1 $, for all $ n\in \N $. By definition, there exists $ \vn \in K_c $ such that $ ||\un - \vn||\leq \frac{1}{n} $, for all $ n \in \N $. 

Now, observe that, setting $ \until = (-[\beta(\un)]) \ast \un $ and $ \vntil = (-[\beta(\vn)]) \ast \vn $, we have $ \until \in L_1 $ and $ \vntil \in K_c $, for all $ n\in \N $. From $ \mathbb{Z} $-invariance, $ ||\until - \vntil||=||\un - \vn|| \leq \frac{1}{n} $.  

Moreover, 
$$
\beta(\vntil)=\beta((-[\beta(\vn)]) \ast \vn)=-[\beta(\vntil)]+\beta(\vntil)
$$
and, from item (iv) of Lemma \ref{l19}, $ \beta(\vntil)\in [0, 1] $. Hence, $ (\vntil)\subset \tilde{K}_c $. Since $ \tilde{K}_c $ is compact in $ \Humeio $, considering a subsequence if necessary, $ \vntil \RA v \in \tilde{K}_c $. From the continuity of $ \beta $, $ \beta(\vntil)\RA \beta(v)\in [0, 1] $. 

On the other hand, 
$$
||\until - v||\leq ||\until - \vntil|| + ||\vntil - v|| \leq \dfrac{1}{n} + o(1) \RA 0,
$$
as $ n \RA + \infty $. Therefore, $ \until \RA v $ in $ \Humeio $. Once more, from $ \beta $ continuity, there exists $ n_0 \in \N $ such that $ \beta(\until)\in [-2, 2] $, for all $ n \geq n_0 $. 

Define $ \wn = (-[\beta(\until)])\ast \vntil $. Since $ \vntil \in K_c $, $ \wn \in K_c $ and $ \beta(\wn)=\beta(\vntil)-[\beta(\until)] $. Consequently,
$$
|\beta(\vntil) - \beta(\until) + (\beta(\until) - [\beta(\until)])| \leq 1 + \dfrac{1}{2}+2 <4 .
$$
Hence, $ \beta(\wn)\in [-4, 4] $, for all $ n \geq n_0 $. So, we conclude that $ \wn \in \tilde{K}_c $, for $ n $ sufficiently large. 

Now, as $ \until \in L_1 $, makes sense to apply $ h_1 $. Thus, 
\begin{align*}
||h_1(\until)-\wn|| & = ||(-[\beta(\until)])\ast \until - (-[\beta(\until)])\ast \vntil|| \\
& = ||(-[\beta(\until)])\ast (\until - \vntil)|| = ||\until - \vntil|| \leq \dfrac{1}{n} .
\end{align*}
Therefore, taking $ n\geq \max\{n_0 , \frac{1}{\overline{\rho}}\} $, we have that $ \wn \in \tilde{K}_c $ and $ ||h_1(\until) -\wn||\leq \overline{\rho} $, which implies that $ h_1(\until)\in U_{\overline{\rho}} $. 

Hence, $ \until \in A_1 $ and, since $ A_1 $ is $ \mathbb{Z} $-invariant, $ \un \in A_1\subset A $, which is a contradiction. So, the claim is valid. 

\noindent \textbf{Claim 4:} $ \gamma(U_{\overline{\rho}})\leq \gamma(\tilde{K}_c) $.

If $ \gamma(\tilde{K}_c)=0 $, then $ \tilde{K}_c= \emptyset $, which implies $ U_{\overline{\rho}}=\emptyset $ and, by definition, $ \gamma(U_{\overline{\rho}})=0 $. Suppose then, $ \gamma(\tilde{K}_c)>0 $. So, we can consider a function $ \tilde{g}:U_{\overline{\rho}}\RA \mathbb{R}^k \{0\} $, given by $ \tilde{g}(u)=g(u) $. Observe that $ \tilde{g} $ is well-defined, since $ U_{\overline{\rho}}\subset \Humeio $ and $ g\neq 0 $ over $ U_{\overline{\rho}} $. Moreover, as $ g $ is odd and continuous, so is $ \tilde{g} $.  Consequently, from genus definition, $ \gamma(U_{\overline{\rho}})\leq k = \gamma(\tilde{K}_c) $, proving the claim.

Finally, since $ \tilde{K}_c , A_{c, \rho}, U_{\overline{\rho}}, A_i, A \in \A $, for $ i=1, 2 $, make sense to apply $ \gamma $ in all of this sets. Therefore, from the upper considerations, Claim 3, and standard properties of genus, for $ \rho \in (0, \rho_c) $, we have
\begin{align*}
\gamma(A_{c, \rho})\leq \gamma(A) & = \gamma\left(A_1 \cup A_2 \right) \\
& \leq \gamma(A_1) + \gamma(A_2) = \gamma(h_{1}^{-1}(U_{\overline{\rho}})) + \gamma(h_{2}^{-1}(U_{\overline{\rho}}))\\
& \leq \sum\limits_{i=1}^{2}\gamma(U_{\overline{\rho}}) \leq  \sum\limits_{i=1}^{2}\gamma(\tilde{K}_c)=2 \gamma(\tilde{K}_c) = 2k < \infty .
\end{align*}
Hence, taking $ \rho_0 = \rho_c $, we have the proposition.
\end{proof}

For the next results we will need the definition and some basic properties of relative genus. So, for convenience of the reader, we will include it here.

\bd\label{d11}
Let $ D, Y \in \A $ with $ D \subset Y $. We say that $ U, V\in \A $ is a covering of $ Y $ relative to $ D $ if it satisfies

\noindent \textbf{(i)} $ Y\subset U \cup V $ and $ D\subset U $;

\noindent \textbf{(ii)} there exists an even continuous (in $ X $) function $ \chi : U \RA D $, such that $ \chi(u)=u $, for all $ u\in D $. 

If $ U, V\in \A $ is a covering of $ Y $ relative to $ D $, then the genus of this covering is $ \gamma(V)=k $. 
\ed 

\bd\label{d12}
Let $ D, Y \in \A $ with $ D \subset Y $. We define the Krasnoselskii's genus of $ Y $ relative to $ D $, denoted by $ \gamma_D(Y) $, as

\noindent \textbf{(i)} There exists a covering for $ Y $ relative to $ D $ and, in this case, $ \gamma_D(Y)=k $, where $ k $ is the lowest genus of this coverings. 

\noindent \textbf{(ii)} If we cannot find any such covering of $ Y $ relative to $ D $, we set $ \gamma_D(Y)=+\infty $.
\ed

In the following, we list some useful properties of relative genus that are needed to guarantee that the results are valid. 

\bl\label{l20}
\textbf{(i)} Let $ D \subset \A $. Then, $ \gamma_D(D)=0 $.

\noindent \textbf{(ii)} Let $ D, Y, Z \in \A $ satisfying $ D\subset Y $ and $ D\subset Z $. If there exists a function $ \varphi : Y \RA Z $, even and continuous (in $ X $), such that $ \varphi(u)=u $, for all $ u\in D $, then $ \gamma_D(Y)\leq \gamma_D(Z) $.

\noindent \textbf{(iii)} Let $ D\subset Y\subset Z \in \A $. Then, $ \gamma_D(Y)\leq \gamma_D(Z) $.

\noindent \textbf{(iv)} Let $ D, Y, Z \in \A $ satisfying $ D\subset Y $. Then, $ \gamma_D(Y\cup Z)\leq \gamma_D(Y)+\gamma(Z) $.
\el
\begin{proof}
For item (i), take $ U=D , V= \emptyset $ and $ \chi = id $, in the definition of relative genus. Proofs for itens (ii) and (iv) can be found, for example, in \cite{[6]}. Finally, item (iii) is an immediate consequence of item (ii).
\end{proof}

For the next results, we define the sets
$$
I^{c}=\{u\in X \ ; \ I(u)\leq c \} \mbox{ \ , for \ } c\in \R \mbox{ \ and \ } D=I^0 ,
$$
and the values
$$
c_k = \inf\{c \geq 0 \ ; \ \gamma_D(I^c)\geq k \} \mbox{ \ , \ } \forall \ n\in \N .
$$

\bo\label{obs14}
\textbf{(1)} Since $ I $ is unbounded from bellow, $ D\neq \emptyset $. 

\noindent \textbf{(2)} Let $ c_1, c_2 \in \R $ with $ c_1 > c_2 $. Then, if $ u\in I^{c_2} $, $ I(u)\leq c_2 < c_1 $, so $ u\in I^{c_2} $. That is, if $ c_1 > c_2 $, then $ I^{c_2}\subset I^{c_1} $.

\noindent \textbf{(3)} If $ c_1 > c_2 \geq 0 $, then $ D\subset I^{c_2}\subset I^{c_1} $. Consequently, $ \gamma_D(I^{c_2})\leq \gamma_D(I^{c_1}) $.

\noindent \textbf{(4)} For $ \varepsilon >0 $, $ \gamma_D(I^{c_k + \varepsilon})\geq k $ and $ \gamma_D(I^{c_k - \varepsilon})<k $, for every $ k \in \N $.
\eo

\bl\label{l112}
We have that
\begin{equation}\label{eq120}
\inf\limits_{u\in X \setminus \{0\}} \sup\limits_{t \in \R} I(tu) = \inf\limits_{u\in X \setminus \{0\}} \sup\limits_{t >0 } I(tu) < + \infty .
\end{equation}
\el
\begin{proof}
First of all, observe that, since $ I $ is even, $\sup\limits_{t \in \R} I(tu) =  \sup\limits_{t >0 } I(tu) $. Moreover, from Lemma \ref{l8}, $ \sup\limits_{t>0} I(tu) < + \infty $ for all $ u \in X\setminus \{0\} $.
\end{proof}

Consider the Nehari's manifold for $ I $, defined by
\begin{equation}\label{VN}
\Nn = \{ u\in X \setminus \{0\} \ ; \ I'(u)(u) = 0 \}.
\end{equation}

\bl\label{l113}
Let $ \Nn $ as in (\ref{VN}). Then, 
$$
\inf\limits_{\Nn} I = \inf\limits_{u \in X \setminus \{0\}}\sup\limits_{t > 0} I(tu) .
$$
\el
\begin{proof}
We start proving the following statement.

\noindent \textbf{Claim 1:} Let $ u\in X \setminus \{0\} $. Then, $ t_u u \in \Nn $ and $ \sup\limits_{t > 0} I(tu) = I(t_u u) $.

Indeed, from Lemma \ref{l110}, there exists an unique $ t_u \in (0, + \infty) $ such that $ \varphi_u '(t_u) =0 $, which, from Corollary \ref{c12}, is a global maximum for $ \varphi_u $, concluding the proof of claim 1.

\noindent \textbf{Claim 2:} Let $ u\in \Nn $. Then, $ \sup\limits_{t > 0} I(tu) = I(u) $.

Since, $ u\in \Nn $ and $ \varphi_u '(t) = I'(tu)(u) $, we have $ 0 = I'(u)(u) = \varphi_u '(1) $. From Lemma \ref{l110}, $ 1 $ is a global maximum for $ \varphi_u $. Consequently, 
$
\sup\limits_{t > 0} I(tu)= \varphi_u (1) = I(u),
$
and the claim follows. 

From Lemma \ref{l110} and Claims 1 and 2, we obtain the result.
\end{proof}

\bl\label{l114}
We have $\inf\limits_{X \setminus \{0\}} \sup\limits_{t\in \R} I(tu) > 0 $. 
\el
\begin{proof}
Suppose that $ \alpha_1 = \inf\limits_{X\setminus\{0\}} \sup\limits_{t\in \R} I(tu) =0 $. Then, there exists a sequence $ (\un) \subset X\setminus \{0\} $ such that $ \sup\limits_{t\in \R} I(t\un) \RA 0 $. From Lemma \ref{l113}, the definition of $ \varphi_{\circ} $ and Lemma \ref{l7} we can reach a contradiction.
\end{proof}

In order to prove the next results, we need to introduce the following sets
$$
\Nn^{+} = \{ u\in X \ ; \ I'(u)(u) > 0 \} \mbox{ \ \ and \ \ } \Nn^{-} = \{ u\in X \ ; \ I'(u)(u) < 0 \}.
$$
Note that $ X = \{0\} \cup \Nn^{+} \cup \Nn \cup \Nn^{-} $ and that this sets are $ 2-2 $ disjoints. Moreover, using the definition of a set's boundary and Lemmas \ref{l7} and \ref{l110}, it is possible to verify that $ \partial \Nn^{-} = \Nn $ and $ \partial \Nn^{+} = \{0\} \cup \Nn $.

\bp\label{p110}
We have $ c_1 = \inf\limits_{N} I = \inf\limits_{u\in X \setminus \{0\}}\sup\limits_{t > 0}I(tu) >0 $.
\ep
\begin{proof}
Observe that, from Lemma \ref{l113}, remains to prove that $ c_1 = \inf\limits_{N} I > 0 $. 

\noindent \textbf{Claim 1:} $ c_1 >0 $. 

Suppose that $ c_1 =0 $. Then, by $ c_1 $ definition, $ \gamma_D(D) \geq 1 $, which contradicts Lemma \ref{l20}-(i). Therefore, $ c_1 >0 $.

\noindent \textbf{Claim 2:} $ c_1 \geq \inf\limits_{\Nn} I $.

Suppose, by contradiction, that $ c_1 < \inf\limits_{\Nn} I  $. Choose $ c\in ( c_1, \inf\limits_{\Nn} I) $. Define the function $ F: I^c \RA X $ by
$$
F(u) = \left\{ \begin{array}{ll}
0 \mbox{ \ , if \ } u\in \{0\} \cup \Nn^{+} \\
\max\{ 1, t_u '\}u \mbox{ \ , if \ } u\in \Nn^{-}
\end{array}
\right. .
$$
Then, we mention some properties of $ F $ that can be easily verify: (i) $ F $ is well-defined, since $ X=\{0\} \cupdot \Nn^{+} \cupdot \Nn \cupdot \Nn^{-} $; (ii) by the definition of $ F $, Lemma \ref{l111}, the continuity of $ \max $ function and $ \partial ( \{0\} \cup \Nn^{+}) = \Nn = \partial \Nn^{-} $, $ F $ is continuous; (iii) from Corollary \ref{c12} and direct computations, one concludes that $ F $ odd; (iv) $ F\big\vert_{D} = id $, it clearly follows from Corollary \ref{c12} and Lemma \ref{l110}.

Consider $ U= I^c $ and $ V = \emptyset $. Since $ I $ is $ C^1 $ and even, $ U $ is closed and  symmetric. That is, $ U\in \A $. Also, $ V = \emptyset \in \A $. Moreover, $ I^c = U \cup V $ and, since $ c\in ( c_1 , \inf\limits_{\Nn} I ) $, $ c> 0 $ and $ D \subset I^c \subset U $. 

From properties (i)-(iv) and definition \ref{d12}, $ U $ and $ V $ is a covering for $ I^c $ relative to $ D $. Then, once $ \gamma(\emptyset) = 0 $, from definition, the genus of this covering is zero. Hence, $ \gamma_D (I^c ) =0 $. But it gives a contradiction, since $ 1 \leq \gamma_D(I^{c_1}) \leq \gamma_D(I^c) = 0 $. Consequently, $ c_1 \geq \inf\limits_{\Nn} I $. 

\noindent \textbf{Claim 3:} $ c_1 \leq \inf\limits_{\Nn} I $.

For $ u_0 \in X\setminus\{0\} $, without loss of generality, we can assume $ ||u_0|| =1 $. 

Set $ d = \sup\limits_{t>0} I(tu_0) $. Lets prove that $ \gamma_D(I^d ) \geq 1 $. Note that, if $ u\in B= \{tu_0 \ ; \ t > 0\} $, then there exists $ t_0 > 0 $ such that $ u= t_0 u_0 $. Then,
$$
I(u) = I(t_0 u_0) \leq \sup\limits_{t>0} I(tu_0) 
$$
and, consequently, $ u\in I^d $. Moreover, $ \inf\limits_{u\in X\setminus\{0\}} \sup\limits_{t>0} I(tu) $ implies $ d>0 $. On the other hand, if $ u\in D $, $ I(u)\leq 0 < d $. So, $ u\in I^d $. Hence, $D \subset B \cup D \subset I^d $. 

Thus, from Lemma \ref{l20}-(iii), $ \gamma_D(B \cup D) \leq \gamma_D(I^d) $. In this point of view, we will work to prove that $ \gamma_D(B \cup D) \geq 1 $. 

Suppose that  $ \gamma_D(B \cup D) =0 $. Then, by definition, $ U= B\cup D $ and $ V = \emptyset $, once only the empty set has null genus. Moreover, there exists a function, continuous and odd, $ \chi : B \cup D \RA D $ such that $ \chi(u)=u $, for all $ u\in D $. 

Define $ g: (0, + \infty) \RA (0, + \infty) $ by $ g(t) = ||\chi(tu_0)|| $. Note that $ g $ is continuous. 

Now, consider $ \rho > 0 $ given by Lemma \ref{l7}. From Corollary \ref{c12}, there exists $ \overline{t} \in (t_{u_0}', + \infty) $, with $ \overline{t} > \rho $, such that $ \varphi_{u_0}(\overline{t}) < 0 $. It is equivalent to $ I(\overline{t} u_0) < 0 $, which implies $ \overline{t}u_0 \in D $. 

Consequently, $ g(\overline{t}u_0) = || \chi(\overline{t}u_0)||= \overline{t}  $. On the other hand, once $ 0\in D $, $ g(0)=0 $. So, since $ g $ is continuous, from the Intermediate Value Theorem, there exists $ t_\rho \in (0, \overline{t}) $ such that $ \rho = g(t_\rho)=||\chi(t_\rho u_0)|| $. But, from Lemma \ref{l7}, $ I(\chi(t_\rho u_0)) \geq m_\rho > 0 $ and we have that $ \chi(t_\rho u_0) \not\in D $, which is a contradiction. 

Therefore, claim 3 is valid and, combined with claim 2, we have the proposition.  
\end{proof}

\bo\label{obs15}
\textbf{(1)} We can provide an equivalent definition to the the function $ F: I^c \RA X $ by
$$
F(u) = \left\{ \begin{array}{ll}
0 \mbox{ \ , if \ } u\in \{0\} \cup \Nn^{+} \\
\sigma(u) u \mbox{ \ , if \ } u\in \Nn^{-}
\end{array}
\right. ,
$$
where $ \sigma : \Nn^{-} \RA  [1, + \infty) $ is given by $ \sigma(u) = \inf\{t \geq 1 \ ; \ \varphi_u(t)=I(tu) \leq 0 \} $. 

\noindent \textbf{(2)} Let $ i, j\in \N $ with $ i>j $. Then, 
$$
\{c\in (0, + \infty) \ ; \ \gamma_D(I^c) \geq i \} \subset \{c\in (0, + \infty) \ ; \ \gamma_D(I^c) \geq j \} . 
$$
Therefore, $ c_j \leq c_i $. 

\eo

In the following, we will consider $ W $ as a $ k $-dimensional subspace of $ X $. As a consequence of $ (f_1 ') $, $ (f_4 ') $, $ V_2(u)\geq 0 $ and $ V_1(u)\leq 2||u||_{X}^{4} $, once can see that 
\begin{equation}\label{l116}
I(u)\RA - \infty \mbox{, as \ }  ||u||_X \RA + \infty \mbox{, and \ } \sup\limits_{u\in W} I(u) < + \infty .
\end{equation}

\bl\label{l117}
There exists $ R>0 $ such that $ \{u\in W \ ; \ ||u|| \geq R \} \subset D $. 
\el
\begin{proof}
Set $ A= \{u\in W \ ; \ ||u|| \geq R \} $. Suppose, by contradiction, that for all $ n\in \N $, there exists $ u_n \in A $ such that $ u_n \not\in D $. Allied with the definition of $ D $, there exists a sequence $ (\un) \subset W $ satisfying $ ||\un||\geq n $ and $ I(\un) > 0 $, for all $ n \in \N $. But it contradicts (\ref{l116}), since for $ n \in \N $ sufficiently large we should have $ I(\un) \leq 0 $.
\end{proof}

For the next results, consider $ \rho > 0  $ given by Lemma \ref{l7} and $ \chi : X \RA X $ a continuous and odd function, such that $ \chi(u)=u $, for all $ u\in D $. Define the sets
$$
\Oo_\chi = \{ u\in W \ ; \ ||\chi(u)|| < \rho \} .
$$ 
For convenience, we set $ \chi $ as a function having these properties, except we say otherwise. 

\bo\label{obs16}
One can easily verify that the sets $ \Oo_\chi $ have the following properties:

\noindent \textbf{(1)} $ \Oo_\chi $ is a neighbourhood of zero in $ W $. 

\noindent \textbf{(2)} $ \Oo_\chi $ is bounded and symmetric.

\noindent \textbf{(3)}  If $ u\in \partial_W \Oo_\chi $, then $ ||\chi(u)||=\rho $.
\eo

\bp\label{p111}
Let $ k\in \N $. Then $ c_k < + \infty $.
\ep
\begin{proof}
Set $ \alpha = \sup\limits_{u\in W} I(u) $. From (\ref{l116}), $ \alpha < + \infty $. Also, by definition of $ I^\alpha $, we see that $ W \subset I^\alpha $. 

Suppose that $ \gamma_D (I^\alpha) < k $. Then, from Definition \ref{d12}, there exists $ U, V \in \A $ and $ \chi : U \RA D $ continuous and odd, such that $ I^\alpha \subset U \cup V $, $ D \subset U $, $ \gamma(V) \leq k-1 $ and $ \chi\big\vert_{D} = id $.

Since $ U\subset X $ is closed, by a corollary of Tietze's Theorem, there exists a unique extension of $ \chi $, that we still denote by $ \chi $, to $ X $, which is continuous and odd. So, the sets $ \Oo_\chi $ are well-defined. 

Now, from Remark \ref{obs16} and a genus property, since $ \dim W = k $, $ \gamma(\partial_W \Oo_\chi ) = k $. 

Let $ u\in \partial_W \Oo_\chi  $. From Remark \ref{obs16}, $ ||\chi(u)||= \rho $. By Lemma \ref{l7}, $ I(\chi(u))\geq m_\rho > 0 $. Thus, $ \chi(u)\not\in D $. Equivalently, $ u\not\in \chi^{-1}(D) $. Once the initial function $ \chi $ satisfies $ \chi(U)\subset D $, we conclude that $ u\not\in U $. Therefore, $ \partial_W \Oo_\chi \cap U = \emptyset $. 

On the other hand, $$ \partial_W \Oo_\chi \subset W \subset I^{\alpha} \subset U \cup V \ \ \Longrightarrow \ \ \partial_W \Oo_\chi \subset V . $$
Then, we have $ k = \gamma(\partial_W \Oo_\chi) \leq \gamma(V) \leq k-1 $, which is a contradiction. 

Therefore, from the definition of values $ c_k $ and (\ref{l116}), $ c_k \leq \alpha < + \infty $. 
\end{proof}

Before we prove that the values $ c_k $ are critical values of $ I $, we will provide a deformation lemma. We start defining the sets $ S = X \setminus A_{c, \rho} $, 
$$
S_\delta = \{ u \in X \ ; \ ||u-v||_X \leq \delta \mbox{ \ , for some \ } v\in S \} 
$$
and 
$$
\tilde{S}_\delta = \{ u \in X \ ; \ ||u-v|| \leq \delta \mbox{ \ , for some \ } v\in S \} ,
$$
for $ c, \rho, \delta \in (0, + \infty ) $. Since the proofs of Lemma \ref{l118} and Corollary \ref{l119} can be done as \cite[Lemma 4.6]{[6]}, we omit it here.

\bl\label{l118}
Let $ c, \rho > 0 $. Then, there exists $ \delta_0 = \delta(c, \rho)>0 $ such that, if $ \delta \in (0, \delta_0) $, then $ ||I'(u)||_{X'}(1+||u||_X) \geq 8 \delta $, for all $ u\in \tilde{S}_{2\delta} $ with $ I(u)\in [c-2\delta^2 , c+ 2 \delta^2] $.
\el

\bc\label{l119}
Let $ c, \rho > 0 $. Then, there exists $ \tilde{\delta}_0 = \tilde{\delta}_0(c, \rho)>0 $ such that, if $ \delta \in (0, \tilde{\delta}_0) $, then $ ||I'(u)||_{X'}(1+||u||_X) \geq 8 \delta $, for all $ u\in S_{2\delta} $ with $ I(u)\in [c-2\delta^2 , c+ 2 \delta^2] $.
\ec

\bl\label{l120}
Let $ c> 0 $. Then, there exists $ \rho_1 = \rho_1 (c) > 0 $ such that, for all $ \rho \in (0, \rho_1)  $, we have

\noindent \textbf{(i)} $ A_{c, \rho} \cap D = \emptyset $ .

\noindent \textbf{(ii)} There exists $ \varepsilon = \varepsilon(c, \rho) >0 $ and a function $ \phi : I^{c+ \varepsilon}\setminus A_{c, \rho} \RA I^{c-\varepsilon} $, continuous and odd, such that $ D \subset I^{c-\varepsilon} $ and $ \phi\big\vert_{D} = id $.
\el
\begin{proof}
\textbf{(i)} Suppose that the statement is false. Then, for all $ n\in \N $, there exists $ \rho_n \in \left(0, \frac{1}{n}\right) $ and $ \un \in A_{c, \rho_n}\cap D $. Note that $ \rho_n \RA 0 $. Thus, by $ A_{c, \rho_n} $ definition, there exists $ \vn \in K_c $ such that $ ||\un - \vn || \leq \rho_n $, for all $ n \in \N $. 

As a consequence, $ I'(\vn) = 0 $ and $ I(\vn)=c $, for all $ n\in \N $. Hence, $ (\vn) $ is a Cerami sequence for $ I $ in level $ c>0 $. Then, by Proposition \ref{p1}, there exists $ (\yn) \subset \mathbb{Z} $ satisfying $ \vntil = \yn \ast \vn \RA v $ in $ X $, where $ v $ is a non-zero critical point of $ I $. Therefore, 
$$
||\until - v|| \leq ||\until - \vntil|| + ||\vntil - v|| = ||\un - \vn||+ ||\vntil - v|| \RA 0 .
$$
Finally, as $ K_c $ is closed, $ v\in K_c $ and, from the $ \mathbb{Z} $-invariance of $ I $ and $ A_{c, \rho_n} $, $ \until = \yn \ast \un  \in A_{c, \rho_n}\cap D  $. So, since $ (\until) \subset D $, $ I(\until) \leq 0 $, for all $ n \in \N $. Consequently, from Lemma \ref{l3}-(iv), 
$$
c = I(v) \leq \liminf I(\until) \leq 0 ,
$$ 
which is a contradiction. 

\noindent \textbf{(ii)} Let $ \delta_0 $ as given by Corollary \ref{l119} and $ \delta\in (0, \delta_0) $, such that $ \delta^2 < \frac{c}{2} $. Take $ \varepsilon = \delta^2 $. Then, from Deformation's Lemma 2.6 of \cite{[def]}, there exists $ \eta : [0, 1]\times X \RA X $, continuous, satisfying

\noindent \textbf{(a)} $ \eta(t, u)=u $, if $ t=0 $ or $ u\not\in I^{-1}([c-2\varepsilon , c+ 2 \varepsilon])\cap S_{2 \delta} $ ;

\noindent \textbf{(b)} $ \eta(1, I^{c+ \varepsilon}\cap (X \setminus A_{c, \rho})) \subset I^{c-\varepsilon} $;

\noindent \textbf{(c)} $ t \mapsto I(\eta(t, u)) $ is non-increasing, for all $ u\in X $.

Moreover, since $ I $ is even, it is possible to modify the proof of the refereed lemma, such as in \cite{[f19]}, to obtain as well

\noindent \textbf{(d)} $ \eta(t, -u) = - \eta(t, u) $, for all $ t \in [0, 1] $ and $ u\in X $.

Define $ \phi : I^{c+ \varepsilon}\setminus A_{c, \rho} \RA I^{c-\varepsilon}  $ by $ \phi(u)=\eta(1, u) $. Note that item (b) is equivalent to $ \eta(1, I^{c+ \varepsilon}\cap \setminus A_{c, \rho}) \subset I^{c-\varepsilon} $, which guarantee that $ \phi $ is well-defined. Also, as $ \eta $ is continuous, $ \phi $ is continuous as well, and from item (d), $ \phi $ is odd. 

Moreover, once $ \varepsilon = \delta^2 < \frac{c}{2} $, $ 0\not\in [ c-2 \varepsilon , c+ 2 \varepsilon] $. Thus, $ D \cap I^{-1}([ c-2 \varepsilon , c+ 2 \varepsilon]) = \emptyset $. Hence, from item (a), if $ u\in D $, $ \phi(u)=\eta(1, u)=u $, which implies that $ \phi\big\vert_{D}=id $. 

Finally, if $ u\in D $, $ I(u)\leq 0< c-\varepsilon $, then $ u\in I^{c-\varepsilon} $. 
\end{proof}

\bp\label{p112}
Let $ k\in \N $. Then, $ c_k $ is a critical value of $ I $.
\ep
\begin{proof}
Suppose, by contradiction, that for all $ u\in X $, with $ I'(u)=0 $, $ I(u)\neq c_k $, for all $ k\in \N $. Let $ \rho > 0 $. From definition of $ K_{c_k} $, we have $ K_{c_k}= \emptyset $ and, consequently, $ A_{c_k, \rho} = \emptyset $, for all $ \rho >0 $. 

From Lemma \ref{l120}, there exists $ \varepsilon=\varepsilon(c_k , \rho) >0 $ and $ \phi: I^{c_k + \varepsilon}\RA I^{c_k - \varepsilon} $, continuous and odd, such that $ \phi\big\vert_{D}=id $. Then, from Lemma \ref{l20}-(ii), $ \gamma_D(I^{c_k + \varepsilon}) \leq \gamma_D(I^{c_k - \varepsilon}) $. 

On the other hand, as $ c_k - \varepsilon < c_k + \varepsilon $, $D\subset I^{c_k - \varepsilon}\subset I^{c_k + \varepsilon}$ and then, by Lemma \ref{l20}-(iii), $ \gamma_D(I^{c_k + \varepsilon}) \geq \gamma_D(I^{c_k - \varepsilon}) $. Consequently, $ \gamma_D(I^{c_k + \varepsilon}) = \gamma_D(I^{c_k - \varepsilon}) $, which is a contradiction, since $ \gamma_D(I^{c_k + \varepsilon})\geq k $ and $ \gamma_D(I^{c_k - \varepsilon})<k $. 

Therefore, $ c_k $ is a critical value of $ I $.
\end{proof}

\bp\label{p113}
We have $ c_k \RA + \infty $, as $ k \RA + \infty $.
\ep
\begin{proof}
Suppose that there exists $ M >0 $, such that $ c_k < M $, for all $ k\in \N $. From Remark \ref{obs15}, $c_k$ is monotonically nondecreasing. Then, there exists $ c > 0 $ such that $ c_k \RA c $.

From Proposition \ref{p19}, there exists $ \rho_0 > 0 $ such that $ \gamma(A_{c, \rho}) < + \infty $, for all $ \rho \in (0, \rho_0) $. Also, by Lemma \ref{l120}-(i), there exists $ \rho_1 > 0 $ such that $ A_{c, \rho}\cap D = \emptyset $, for all $ \rho \in (0, \rho_1) $. Set $ \rho_2 = \min\{\rho_0 , \rho_1\} $ and let $ \rho \in (0, \rho_2 )$.

From Lemma \ref{l120}-(ii), there exists $ \varepsilon >0 $ and $ \phi : I^{c+ \varepsilon}\setminus A_{c, \rho} \RA I^{c-\varepsilon} $, continuous and odd, such that $ \phi\big\vert_{D}=id $ and $ D\subset I^{c- \varepsilon} $. Moreover, once $ c_k \RA c $ monotonically nondecreasing, there exists $ k_0 \in \N $, such that $ c_{k_0}\geq c-\varepsilon $.  Then, by Lemma \ref{l20}-(ii), $$ \gamma_D(I^{c+ \varepsilon}\setminus A_{c, \rho}) \leq \gamma_D(I^{c-\varepsilon}) < k_0 < + \infty . $$
Consequently, from Lemma \ref{l20}-(iv), 
$$
\gamma_D (I^{c+ \varepsilon}) = \gamma_D((I^{c+ \varepsilon}\setminus A_{c, \rho})\cup A_{c, \rho}) \leq \gamma_D (I^{c+ \varepsilon}\setminus A_{c, \rho}) + \gamma(A_{c, \rho}) < + \infty .
$$
Once more, as $ c_k \RA c $ monotonically nondecreasing, $ c+ \varepsilon > c \leq c_k $, for all $ k \in \N $. Also, by definition of values $ c_k $, $ \gamma_D (I^{c+ \varepsilon}) \geq k $, for all $ k \in \N $, that is $ \gamma_D (I^{c+ \varepsilon}) \RA + \infty $, as $ k \RA + \infty $, which leads to a contradiction. 

Therefore, $ c_k \RA + \infty $, as $ k \RA + \infty $.
\end{proof}

\begin{proof}[Proof of Theorem \ref{t2}] 
From Proposition \ref{p113}, we can extract a subsequence of $ (c_k) $ such that $ c_k \RA + \infty $ monotonically increasing. Then, from Proposition \ref{p112}, there exists $ \uk \in X $ satisfying $ I(\uk)=c_k $ and $ I'(\uk)=0 $, for all $ k\in \N $. Since $ (c_k) $ is monotone, $ c_i \neq c_j $, when $ i\neq j $, and $ c_k > c_1 >0 $, for all $ k\geq 2 $. Then, we conclude that the functions $ u_k $ are distinct and that $ \uk \neq 0 $, for all $ n\in \N $. Also, as $ I $ is even, the same holds for $ -\uk $ and $ I(\pm \uk) \RA + \infty $. 
\end{proof}

\end{document}